\documentclass[12pt]{article}
\usepackage{amssymb,amsmath}
\usepackage{epsfig}

\newtheorem{Theorem}{Theorem}[section]
\newtheorem{Proposition}[Theorem]{Proposition}

\newtheorem{Lemma}[Theorem]{Lemma}

\newenvironment{proof}{\par\medskip\noindent{\it Proof.}~}{\hfill $\square$}

\def\<{\langle}
\def\>{\rangle}
\def\Z{\mathbb{Z}}

\title{Magnus subgroups of one-relator surface groups}

\author{James Howie and Muhammad Sarwar Saeed\\
Department of Mathematics and\\ Maxwell Institute of Mathematical Sciences\\
Heriot-Watt University\\ Edinburgh EH14 4AS\\ UK}

\begin{document}

\maketitle

\begin{abstract}
A one-relator surface group is the quotient of an orientable
surface group by the normal closure of a single relator.  A Magnus
subgroup is the fundamental group of a suitable incompressible
sub-surface.  A number of results are proved about the
intersections of such subgroups and their conjugates, analogous to
results of Bagherzadeh, Brodski\u\i\, and Collins in classical
one-relator group theory.
\end{abstract}

\section{Introduction}

Recall the Freiheitssatz of Magnus  \cite{Mg30, Mg31} for one-relator groups:

\begin{Theorem} {\rm [The Freiheitssatz]}
Let $G=\langle X:R\rangle$ be a one-relator group where $R$ is
cyclically reduced. If $Y$ is a subset of $X$ which omits a
generator occurring in $R$, then the subgroup $M_Y$ generated by
$Y$ is freely generated by $Y$.
\end{Theorem}

Subgroups of a one-relator group of the form $M_Y$ as in the
Freiheitssatz are called {\em Magnus subgroups}. In \cite{New68},
Newman proved that the Magnus subgroups of a one-relator group
with torsion are malnormal, that is, if $M$ is a Magnus subgroup
and $g\notin M$, then $M\cap gMg^{-1}$ is trivial. Bagherzadeh
\cite{Ba76} generalized Newman's result in 1976 to ordinary
one-relator groups and proved that Magnus subgroups of one-relator
groups are cyclonormal. He proved the following

\begin{Theorem}\label{b}
Let $M$ be a Magnus subgroup of a one-relator group $G=\langle
X:R\rangle$. Then $M$ is cyclonormal in $G$, that is, if $g\notin
M$, then $M\cap gMg^{-1}$ is cyclic.
\end{Theorem}

Collins \cite{Col04,Col05} proved the following results about the
intersection of Magnus subgroups of a one-relator group $G$.

\begin{Theorem}\label{col1}
Let $M_Y$ and $M_Z$ be Magnus subgroups of a one-relator group
$G=\langle X:R\rangle $ generated by subsets $Y,Z\subset X$. Then
$$M_Y\cap_G M_Z=M_{Y\cap Z}\ast I,$$ where $I$ is a free group of
rank $0$ or $1$.
\end{Theorem}

\begin{Theorem}\label{col2}
Let $M_Y$ and $M_Z$ be Magnus subgroups of a one-relator group $G$
as in Theorem \ref{col1}. For any $g\in G$, either $M_Y\cap_G
gM_Zg^{-1}$ is cyclic (possibly trivial) or $g\in M_YM_Z$.
\end{Theorem}

Here we use the notation $A\cap_G B$ to denote the intersection of
two subgroups $A,B$ in the group $G$, to distinguish it from the
intersection in any other group containing them both.  For
example, in Theorem \ref{col1}, if $F$ is the free group on $X$,
then $M_Y\cap_F M_Z=M_{Y\cap Z}$; the Theorem tells us that this
may differ from $M_Y\cap_G M_Z$.

When these two intersections do differ, in other words when $I$
has rank 1 in Theorem \ref{col1}, we say that the two Magnus
subgroups involved have {\it exceptional intersection}. The first
author \cite[Theorem E]{How05} has shown that it is
algorithmically decidable whether a given pair of Magnus subgroups
in a given one-relator group has exceptional intersection.

A {\em one-relator surface group} is the quotient of the
fundamental group of an orientable surface (possibly noncompact,
or with boundary) by the normal closure of a single element. These
groups were introduced in 1990 by Hempel \cite{Hem90}, and have
subsequently been studied by Bogopolski and Sviridov \cite{Bog,BS}
and by the first author \cite{How}.

In this paper we generalize Theorems \ref{b}, \ref{col1} and
\ref{col2}, as well as \cite[Theorem E]{How05} from one-relator
groups to one-relator surface groups.  With an appropriate
definition of Magnus subgroup, we prove the following.

\medskip\noindent{\bf Theorem \ref{s}}
{\it Let $G$ be a one-relator surface group, and let $M$ be a
Magnus subgroup of $G$. Then $M$ is cyclonormal, that is, for any
$g\in G\smallsetminus M$, $M\cap gMg^{-1}$ is cyclic.}

\medskip\noindent{\bf Theorem \ref{ss}}
{\it The intersection $M_1\cap_GM_2$ of two compatible Magnus
subgroups $M_1$ and $M_2$ of the one-relator surface group $G$ is
the free product of $(M_1\cap_{\Sigma}M_2)$ with a cyclic group.
That is, $M_1\cap_GM_2=(M_1\cap_{\Sigma}M_2)\ast C$.}

\medskip
(See \S \ref{prelim} for the definition of {\em compatible} Magnus subgroups.)

\medskip\noindent{\bf Theorem \ref{alg}}
{\it There is an algorithm which will decide, given a one-relator surface
group and a pair of compatible Magnus subgroups, whether or not the
intersection is exceptional (that is, $C\ne\{1\}$ in Theorem \ref{ss})
and if so will give a generator for $C$.}

\medskip\noindent{\bf Theorem \ref{sss}}
{\it Let $G$ be a one-relator surface group and let $M_1$ and
$M_2$ be two Magnus subgroups of $G$. Let $g\in G$, then
$M_1\cap_G gM_2g^{-1}$ is cyclic unless $g\in M_1M_2$.}

Theorems \ref{s}, \ref{ss} and \ref{sss} appeared in the second
author's thesis \cite{MSS}. We are grateful to the thesis
exmainers, Andrew Duncan and Nick Gilbert, for useful comments.

\medskip
In \S\ref{prelim} below we define our notion of Magnus subgroup for
one-relator surface groups and present some useful preliminary results.
Theorem \ref{s} is proved in \S\ref{cyclo}, Theorems \ref{ss} and \ref{alg}
in \S\ref{except}, and finally Theorem \ref{sss} in \S\ref{conj}.

\section{Preliminaries}\label{prelim}

In order to formulate appropriate generalizations of theorems
about Magnus subgroups of one-relator groups, we first need to
choose a suitable definition of Magnus subgroup for a one-relator
surface group. A minimum requirement for a Magnus subgroup is that
it should satisfy an appropriate version of the Freiheitssatz for
one-relator surface groups - which turns out to be a somewhat
delicate question (see \cite{How}).  For the purposes of
exposition in the present paper we shall restrict our definition
of Magnus subgroup to a case where we know that a Freiheitssatz
holds.

Suppose that $S$ is a surface and $\alpha$ is an essential separating
simple closed curve on $S$.  Then the surface group $\Sigma=\pi_1(S)$
splits along $\alpha$ as a free product with amalgamation:
$$\Sigma=\pi_1(S)\cong\pi_1(S_1)\ast_A\pi_1(S_2),$$
where $S_1$ and $S_2$ are the two components of $S$ cut along
$\alpha$, and $A$ is the cyclic subgroup generated by $\alpha$.

Provided that $R\in\pi_1(S)$ is not conjugate into one of the
factors $\pi_1(S_i)$, we say that $\pi_1(S_1)$ and $\pi_1(S_2)$
are {\em Magnus subgroups} in the one-relator surface group
$\pi_1(S)/\<\<R\>\>$.  Note that a Magnus subgroup is generated by
a subset of some standard generating set for the surface group
$\Sigma=\pi_1(S)$ -- for example
$$\<a_1,b_1,\dots,a_{\ell},b_{\ell}\>\subset
\Sigma=\<a_1,b_1,\dots,a_k,b_k~:~[a_1,b_1]\cdots [a_k,b_k]=1\>.$$

\begin{Theorem}{\rm (Freiheitssatz for one-relator surface groups \cite[Proposition 3.10]{How})}\label{fhs}
If $M=\pi_1(S_1)$ is a Magnus subgroup in a one-relator surface group
$G$, then the inclusion map $M\to G$ is injective,
\end{Theorem}

The separating curve $\alpha$ in the definition is determined by
the Magnus subgroup only up to isotopy. We shall also refer to a
pair of Magnus subgroups $M_1$ and $M_2$ as {\em compatible} if
the corresponding separating curves on $S$ can be chosen to be
disjoint.  In terms of generators, there exists a standard
generating set such that both $M_1,M_2$ are generated by subsets
of the chosen generating set for $\Sigma=\pi_1(S)$.

\medskip
\noindent
{\bf Remark} Let $$G=\Sigma/\<\<R\>\>=\<a_1,b_1,\ldots ,a_k,b_k~:~[a_1,b_1]\cdots [a_k,b_k]=R=1\>$$
 be a one-relator surface group, and $L=\{a_1,b_1,\ldots ,a_{k-1},b_{k-1},b_k\}$ a
proper subset of the generating set of $G$. Then $L$ generates a
subgroup $M$ of $\Sigma=\pi_1(S)$ corresponding to the complement
of a nonseparating simple closed curve in $S$. In \cite{How} it is
shown that the Freiheitssatz does not in general hold for such
subgroups: the natural map $M\to G$ is not always injective. For
this reason, we have excluded such subgroups from our definition
of {\em Magnus subgroup}.

Moreover, it turns out that the results of this paper do not necessarily
extend to groups of this form, even in situations where $M\to G$
is injective.  We shall give an example in \S \ref{cyclo} to illustrate this.

In \S \ref{cyclo} below we will employ an idea first used by
Hempel \cite[Lemma 2.1, Theorem 2.2]{Hem90} (see also
\cite[Proposition 2.1]{How}) to express a one-relator surface
group as an HNN extension of a one-relator group. Here the
notation $\<\alpha,\beta\>$ denotes the algebraic intersection
number of a pair of curves $\alpha,\beta$ on the surface $S$.

\begin{Proposition}\label{trick}
Let $S$ be a closed, connected, oriented surface of genus at least $2$,
and let $\alpha$ be a closed curve in $S$.  Then
\begin{enumerate}
\item There is a non-separating simple closed curve $\beta$ in $S$
such that $\<\alpha,\beta\>=0$. \item For any such $\beta$, there
are connected surfaces $F,F_0,F_1$ and a closed curve $\alpha'$ in
$F$, such that
\begin{enumerate}
\item[(a)] $F_0\cong F_1$, $F_0\subset F$ and $F_1\subset F$;
\item[(b)] $\pi_1(F_0)\to\pi_1(F)/\<\<\alpha'\>\>$ and $\pi_1(F_1)\to\pi_1(F)/\<\<\alpha'\>\>$ are injective;
\item[(c)] $\pi_1(S)$ (resp. $\pi_1(S)/\<\<\alpha\>\>$) is an HNN-extension of
$\pi_1(F)$ (resp. $\pi_1(F)/\<\<\alpha'\>\>$) with associated subgroups $\pi_1(F_0)$ and $\pi_1(F_1)$;
\item[(d)] Each of $\partial F$,  $\partial F_0$ and $\partial F_1$
consists of two circles, each of which represents (a conjugate of)
$\beta\in\pi_1(S)$.
\end{enumerate}
\end{enumerate}
\end{Proposition}

In \S \ref{except} and \S\ref{conj} we will use a slight variation
of this idea, which we will describe in the course of the proof of Theorem \ref{ss}.

We shall also make extensive use of the fact that there is a lot
of freedom in the choice of the curve $\beta$. In particular, if
$S_0\subset S$ is a punctured torus, then the restriction of the
algebraic intersection map $\<\alpha,-\>$ to $S_0$ gives a
homomorphism $\Z^2\cong H_1(S_0)\to\Z$ with nonzero kernel; we may
choose a simple closed curve $\beta\subset S_0$ to represent a
nonzero element of the kernel, and such a curve is automatically
nonseparating.

Lemma \ref{sl} below is an algebraic translation of this
observation, applied to the case of the closed orientable surface
$S$ of genus $g$, with a standard generating set
$\{a_1,b_1,\dots,a_k,b_k\}$ for $\pi_1(S)$, where $\pi_1(S_0)$ is
generated by $\{a_k,b_k\}$.

If $R$ is an element of a free group $F$ with basis $X$, and $x\in X$,
we denote by $\sigma(R,x)$ the exponent-sum of $x$ in $R$, in other
words the image of $R$ under the homomorphism $F\to\Z$ defined by
$x\mapsto 1$, $X\smallsetminus\{x\}\mapsto 0$.

\begin{Lemma}\label{sl}
Let $R$ be an element of the free group $\langle a_1,b_1,\ldots
,a_k,b_k\rangle $. Then there exists a basis $\{a_k',b_k'\}$ of
the free group $\langle a_k,b_k\rangle $, such that
\begin{description}
\item (i) $[a_k',b_k']=[a_k,b_k]$; and \item (ii) as a reduced
word in $\{a_1,b_1,\ldots ,a_{k-1},b_{k-1},a_k',b_k'\}$, $R$ has
exponent sum zero in $a_k'$.
\end{description}
\end{Lemma}

\section{Magnus subgroups are cyclonormal}\label{cyclo}

\begin{Theorem}\label{s}
Let $G$ be a one-relator surface group, and let $M$ be a Magnus
subgroup of $G$. Then $M$ is cyclonormal, that is, for any $g\in
G\smallsetminus M$, $M\cap gMg^{-1}$ is cyclic.
\end{Theorem}

\begin{proof}
Without loss of generality, we may assume that
$$G=\<a_1, b_1,\ldots , a_k, b_k~:~[a_1,b_1]\cdots [a_k,b_k]=R=1\rangle$$
and that $$M=\langle a_1, b_1,\ldots , a_{k-1}, b_{k-1}\rangle.$$
Let $g\notin M$ be an element of $G$.

By Lemma \ref{sl}, we may assume that $$\sigma(R,a_k)=0,$$ that
is, the exponent sum of $a_k$ in $R$ is zero.

Let $\beta$ denote the simple-closed curve on $S$ representing
$b_k$, such that $\sigma(-,a_k)=\<-,\beta\>:\pi_1(S)\to\Z$.  Then
Hempel's trick (Proposition \ref{trick}) with this choice of
$\beta$ expresses $G$ as an HNN-extension $$G=\langle H, a_k\mid
a_kXa_k^{-1}=Y\rangle $$ of a one-relator group $H$, in such a way
that $X$ and $Y$ are isomorphic Magnus subgroups of $H$, with $M$
a free factor of $X$ and such that, in fact, $H$ is a one-relator
product of $M$ and $Y$.

The Bass-Serre tree for this HNN-extension has vertex-stabilizers
the conjugates of $H$ and edge-stabilizers the conjugates of $X$
(or the conjugates of $Y$). Let $T$ be the Bass-Serre tree for
this HNN-extension and suppose $g\in G\smallsetminus H$. Then $G$
acts on $T$ and there exists a vertex $v$ such that $$H=Stab(v)$$
$$gHg^{-1}=Stab(g(v)).$$ Moreover, $X$ and $Y$ are the stabilisers
of two edges of $T$, which have $v$ as source and target
respectively.

\begin{figure}[htbp]
\begin{center}
\epsfxsize=10cm \epsfysize=2cm \epsfbox{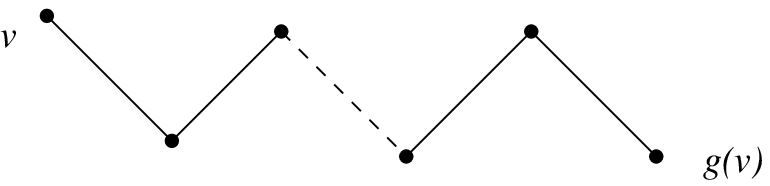}
\end{center}
\caption{}\label{12}
\end{figure}

Now $M\subset H$ stabilizes $v$ and $gMg^{-1}\subset gHg^{-1}$
stabilizes $g(v)$ so that $M\cap gMg^{-1}$ stabilizes both $v$ and
$g(v)$ and hence stabilizes the path $P$ in $T$ from $v$ to
$g(v)$. Here three different cases arise.

\medskip\noindent
{\bf Case 1.} If $g(v)=v$, then $g\in Stab(v)=H$ and the result
follows from Bagherzadeh's Theorem \ref{b}.

\medskip\noindent
{\bf Case 2.} If the path $P$ is not coherently oriented, then
there is an intermediate vertex $u=g'(v)$ of $P$ that is either
the source of each incident edge of $P$ or the target of each
incident edge of $P$.  We treat the latter case (Figure \ref{13});
the former is entirely analogous.

\begin{figure}[htbp]
\begin{center}
\epsfxsize=12cm \epsfysize=2.5cm \epsfbox{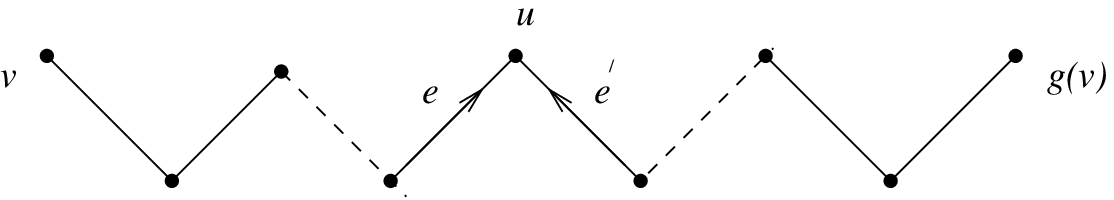}
\end{center}
\caption{}\label{13}
\end{figure}

If $e,e'$ are the edges of $P$ incident at $u$,
then
$$Stab(e)=g'Xg'^{-1}$$
$$Stab(e')=(g'h)X(g'h)^{-1}$$
for some $h\in H$. Now $$M\cap gMg^{-1}\subseteq Stab(v)\cap
Stab(g(v))\subseteq Stab(e)\cap Stab(e').$$But $$Stab(e)\cap
Stab(e')=g'Xg'^{-1}\cap g'hXh^{-1}g'^{-1}=g'(X\cap
hXh^{-1})g'^{-1}.$$ Therefore $$M\cap gMg^{-1}\subseteq g'(X\cap
hXh^{-1})g'^{-1}$$is cyclic by Bagherzadeh's Theorem \ref{b}.

\medskip\noindent
{\bf Case 3.} If the path $P$ is coherently oriented, then we will
assume that the orientation is from $g(v)$ to $v$ -- see Figure
\ref{17}. (The argument for the opposite orientation is
analogous.)

\begin{figure}[htbp]
\begin{center}
\epsfxsize=10cm \epsfysize=2cm \epsfbox{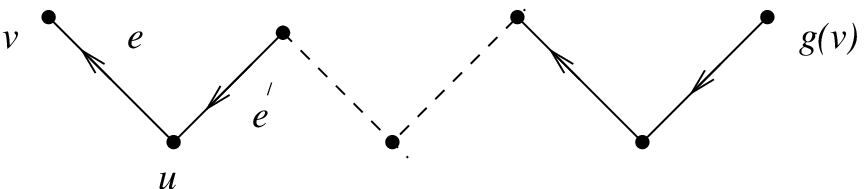}
\end{center}
\caption{}\label{17}
\end{figure}

The edge $e$ of $P$ incident at $v$ has target $v$ and so has
stabilizer $hYh^{-1}$ for some $h\in H$.

Now $$M\cap gMg^{-1}\subseteq M\cap hYh^{-1},$$ and $H$ is a
one-relator product of the free groups $M$ and $Y$, so $M\cap
hYh^{-1}$ is cyclic by a result of Brodski\u\i\, \cite[Teorema
6({\sc b})]{Bi81}.

\medskip
In all cases we have shown that $M\cap gMg^{-1}$ is cyclic. Hence $M$ is cyclonormal in $G$.
\end{proof}

\medskip
Below we give an example to show that Theorem \ref{s} does not extend to
a subgroup $M$ of $G$ generated by $2k-1$ of the $2k$ generators of $G$,
even in cases where $M$ is free on those generators.

\medskip\noindent{\bf Example}
Let $$G=\<
a_1,b_1,a_2,b_2~:~[a_1,b_1][a_2,b_2]=b_1^4a_2^{-1}b_1^3a_2b_1^2a_2^{-1}b_1^3a_2=1\>.$$
Then the second relator $R\equiv
b_1^4a_2b_1^3a_2^{-1}b_1^2a_2b_1^3a_2^{-1}$ has exponent-sum $0$
in $a_2$, so the Freiheitssatz for one-relator surface groups
\cite[Proposition 3.10]{How} implies that $M=\<a_1,b_1,b_2\>$
embeds in $G$ via the natural map. On the other hand, Collins
\cite{Col04} shows that $R=1\Rightarrow b_1^6=a_2^{-1}b_1^6a_2$.
Note also that $[a_1,b_1][a_2,b_2]=1\Rightarrow
a_2^{-1}b_2a_2=b_2[a_1,b_1]$, so that the nonabelian free subgroup
$\<b_1^6,b_2[a_1,b_1]\>$ of $M$ is identified in $G$ with the
subgroup $\<a_2^{-1}b_1^6a_2,a_2^{-1}b_2a_2\>$ of $a_2Ma_2^{-1}$.
Hence $M\cap_G a_2Ma_2^{-1}$ is not cyclic, so $M$ is not
cyclonormal in $G$.

\section{Intersections of Magnus subgroups in one-relator surface groups}\label{except}

Our aim in this section is to prove the analogue of the Theorem
\ref{col1} of Collins for one-relator surface groups.

\begin{Theorem}\label{ss}
The intersection $M_1\cap_GM_2$ of two compatible Magnus subgroups
$M_1$ and $M_2$ of the one-relator surface group $G$ is the free
product of $(M_1\cap_{\Sigma}M_2)$ with a cyclic group. That is,
$$M_1\cap_GM_2=(M_1\cap_{\Sigma}M_2)\ast C.$$
\end{Theorem}

\begin{proof}
Without loss of generality, we may suppose that
$$G=\langle a_1, b_1,\ldots , a_k, b_k~:~[a_1, b_1]\cdots[a_k, b_k]=R=1\rangle,$$
$$M_1= \langle a_1, b_1,\ldots , a_{k-1},b_{k-1}\rangle ,$$
$$M_2=\langle a_2, b_2,\ldots , a_k,b_k\rangle.$$

By Lemma \ref{sl}, we may assume that $a_1, a_k$ appear in $R$
with exponent sum zero, that is, $$\sigma (R,a_1)=0=\sigma
(R,a_k).$$

Note that $$\Sigma =M_1\ast _{M_0}M_2$$ where $M_0=M_1\cap_{\Sigma }M_2$.

\bigskip
By definition of Magnus subgroup, $R$ is not conjugate to an
element of $M_1$ or of $M_2$.  Hence we may assume that $R\in
\Sigma=M_1\ast_{M_0}M_2$ is a cyclically reduced word of length
greater than 1 (with respect to the amalgamated free product
length function).

We apply an amended form of Hempel's trick as follows.  The kernel
$K$ of $\sigma(-,a_k):\Sigma\to\Z$ has an induced graph-of-groups
decomposition as an infinite amalgamated free product of
$\tilde{M_2}=K\cap M_2$ and the groups $a_k^nM_1a_k^{-n}$ for
$n\in\Z$, amalgamating the copy of $a_k^nM_0a_k^{-n}$ in
$a_k^nM_1a_k^{-n}$ with that in $\tilde{M_2}$. Choose a conjugate
$\tilde{R}$ of $R$ that belongs to the subgroup $K_0$ of $K$
generated by $\tilde{M_2}$ and $a_k^nM_1a_k^{-n}$ for $0\le n\le
m$, and assume that all the choices have been made to minimize
$m$.   Then $G$ is an HNN-extension of the one-relator group
$\tilde{G}=K_0/\<\<\tilde{R}\>\>$, with stable letter $a_k$ and
associated subgroups $K_1=K_0\cap a_kK_0a_k^{-1}$, $K_2=K_0\cap
a_k^{-1}K_0a_k$.

\medskip
Clearly $M_1\cap_G M_2\subset M_1\cap_G \tilde{M_2}$. Note  also
that $M_1\subset K_0$.  If $m>0$ in the above construction, then
the join of $M_1$ and $\tilde{M_2}$ in $K_0$ is a Magnus subgroup
of the one-relator group $\tilde{G}$, from which it follows that
$$M_1\cap_G M_2=M_1\cap_{K_0}\tilde{M_2}=M_0.$$

\medskip
Hence we are reduced to the situation where $m=0$ in the HNN
construction. Now $M_0$ is a free factor of $\tilde{M_2}$, so we
can write $\tilde{M_2}=M_0*F$ for some free group $F$, and then we
also have $K_0=M_1*F$.

We now essentially repeat the above argument, with $a_1$ replacing
$a_k$. Specifically, let $N$ be the kernel of
$\sigma(-,a_1):K_0\to\Z$. Then $N$ is the free product of
$\tilde{M_1}=M_1\cap N$ and the groups $F_n:=a_1^nFa_1^{-n}$ for
$n\in\Z$.   Choosing a suitable conjugate $\hat{R}$ of
$\tilde{R}$, we may assume that
$$\hat{R}\in\tilde{M_1}*F_0*F_1*\cdots F_p$$ with all choices made
to minimize $p$. Then $\tilde{G}$ is an HNN-extension of the
one-relator group $\hat{G}=(\tilde{M_1}*F_0*F_1*\cdots
F_p)/\<\<\hat{R}\>\>$ with stable letter $a_1$.

Arguing as before, $M_1\cap_G\tilde{M_2}\subset\tilde{M_1}\cap_G\tilde{M_2}$.
Moreover, $\tilde{M_2}\subset F_0$.  If $p>0$, then the join of
$\tilde{M_1}$ and $F_0$ in $N$ is a Magnus subgroup of $\hat{G}$,
from which it follows that $$\tilde{M_1}\cap_G\tilde{M_2}=\tilde{M_1}\cap_N\tilde{M_2}=M_0.$$

Hence we are reduced to the case where $p=0$.

Now $$\tilde{M_1}=M_0\ast L$$ where  $L$ is a free group.
Also $$\tilde{M_2}=M_0\ast F_0,$$
and
$$\hat G=(M_0\ast F_0\ast L)/\<\<\hat{R}\>\>.$$
Therefore $$M_1\cap_GM_2=\tilde{M_1}\cap \tilde{M_2}=(M_0\ast
F_0)\cap_{\hat{G}} (M_0\ast L).$$ Since $\hat G$ is a one-relator
group, Collins' Theorem  \ref{col1} applies, and so
$$M_1\cap_GM_2=M_0\ast C$$ with $C$ cyclic, as required.
\end{proof}

\medskip
The proof of Theorem \ref{ss}
shows that Magnus subgroups can have exceptional
intersection only in very restricted circumstances - where
$R\in M_0*F_0*L$ in the notation of the proof.  Moreover, in that
case it is equivalent to a pair of Magnus subgroups in a one-relator
group having exceptional intersection.  We can use this to generate
examples of exceptional intersections of Magnus subgroups in
one-relator surface groups.

\medskip\noindent{\bf Example}
Let $G$ be the one-relator surface group
$$\< a_1,b_1,a_2,b_2 ~:~ [a_1,b_1][a_2,b_2]=a_1^{-2}b_1^4a_1^2a_2^{-2}b_2^{-3}
a_2^2a_1^{-2}b_1^2a_1^2a_2^{-2}b_2^{-3}a_2^2=1\>.$$
If $x=a_1^{-2}b_1a_1^2$ and $y=a_2^{-2}b_2a_2$, then the second relation
is $x^4y^{-3}x^2y^{-3}=1$.  Collins \cite{Col04} shows that $x^6=y^6$
is a consequence of that relation.
If $M_1=\<a_1,b_1\>$ and $M_2=\<a_2,b_2\>$, then $x^6\in M_1$ and
$y^6\in M_2$.   Hence $M_1$ and $M_2$ have exceptional
intersection in $G$.

\medskip The strong restrictions on exceptional intersection that
arise in the proof of Theorem \ref{ss} also give rise to a proof of
Theorem \ref{alg}, which we sketch below

\begin{Theorem}\label{alg}
There is an algorithm which will decide, given a one-relator surface
group $G$ and two Magnus subgroups $M_1,M_2$, whether or not $M_1$
and $M_2$ have exceptional intersection in $G$.  If the intersection
is exceptional, the algorithm will provide a generator for the free
factor $C$ in the statement of Theorem \ref{ss}.
\end{Theorem}

\medskip\noindent{\em Sketch proof.}
The theorem is proved by noting that each step in the proof of Theorem \ref{ss}
can be carried out algorithmically.

We may assume that the one-relator surface group has the form
$$\<a_1,b_1,\dots,a_k,b_k~:~[a_1,b_1]\cdots[a_k,b_k]=R=1\>,$$
where $R$ is a word in the generators.

The first step is a basis-change in the free group $\<a_1,b_1\>$
to allow us to assume that $\sigma(R,a_1)=0$.  The euclidean
algorithm transforms the vector $(\sigma(R,a_1),\sigma(R,b_1))\in\Z^2$
to a vector of the form $(0,\ell)$ using integer elementary column
operations, which can be lifted to Nielsen operations on $\<a_1,b_1\>$
in the standard way.  Thus the basis-change operation of Lemma \ref{sl}
can be performed algorithmically, and so we may assume without further
ado that $\sigma(R,a_1)=0$, and similarly that $\sigma(R,a_k)=0$.

The rewrites $R\to\tilde{R}\to\hat{R}$ in the proof of Theorem \ref{ss}
are entirely mechanical processes, as is the choice of a suitable
cyclic conjugate in each case.  Thus the non-negative integers
$m,p$ occurring in the proof can be algorithmically computed.
Should either be strictly positive, then we can stop, declaring the
intersection to be non-exceptional.

Hence we may assume that $m=p=0$, so that (up to conjugation),
$R\in M_0*F_0*L$ in the notation of the proof of Theorem \ref{ss}.
Now $F_0$ and $L$ are free groups of infinite rank, so in order to
handle this situation algorithmically we must replace them by appropriate
finite rank free groups.  In practice, one can algorithmically generate
finite sets $B_1,B_2$ that are subsets of bases of $F_0,L$ respectively,
and such that $R$ can be expressed (up to conjugacy) as a word in
$M_0*\<B_1\>*\<B_2\>$.

Now apply the algorithm of \cite[Theorem E]{How05} to the one-relator group
$(M_0*\<B_1\>*\<B_2\>)/\<\<R\>\>$ to decide whether or not the intersection
is exceptional.  If so, the algorithm provides a generator $\gamma$ for the
exceptional free factor, in terms of our chosen basis for $M_0$
together with $B_1\cup B_2$.  Finally, we translate $\gamma$ into a
word in the original generators $a_1,b_1,\dots,a_k,b_k$ of $G$ to complete
the algorithm.

\section{Intersections of conjugates of Magnus subgroups of one-relator surface groups}\label{conj}

 In this section we prove the analogue of Theorem \ref{col2}.

\begin{Theorem}\label{sss}
Let $G$ be a one-relator surface group and let $M_1$ and $M_2$ be
two compatible Magnus subgroups of $G$. Let $g\in G$. Then $M_1\cap_G
gM_2g^{-1}$ is cyclic unless $g\in M_1M_2$.
\end{Theorem}

\begin{proof}
As in the proof of Theorem \ref{ss}, we assume that
$$G=\langle a_1, b_1,\ldots , a_k, b_k:[a_1, b_1]\cdots[a_k, b_k]=R=1\rangle,$$
$M_1=\<a_1, b_1,\ldots , a_{k-1},b_{k-1}\>$ and $M_2=\<a_2,
b_2,\ldots , a_k,b_k\>$. We also assume, by virtue of Lemma
\ref{sl}, that $\sigma(R,a_1)=\sigma(R,a_k)=0$.

Let $g\in G$.  Note that for any $m,n\in \mathbb{Z}$ we may
replace $g$ by $g'=a_1^mga_k^n$, since $M_1\cap
g'M_2(g')^{-1}=M_1\cap gM_2g^{-1}$. Hence we may assume that
$\sigma(g,a_1)=0=\sigma(g,a_k)$.

With the same notation as in the proof of Theorem \ref{ss}, we
express $G$ as an HNN extension of a one-relator group
$\tilde{G}=K_0/\<\<\tilde{R}\>\>$, with stable letter $a_k$ and
associated subgroups $K_1$, $K_2$, where  $K_0$ is generated by
$\tilde{M_2}=M_2\cap\mathrm{Ker}(\sigma(-,a_k))$ together with
$a_k^nM_1a_k^{-n}$ for $0\le n\le m$, for some $m\ge 0$. In
particular $M_1\subset K_0$.

Note that, since $M_1\subset\mathrm{Ker}(\sigma(-,a_k))$, we have
$$M_1\cap_G gM_2g^{-1}=M_1\cap_G g\tilde{M_2}g^{-1}\subset \tilde{G}\cap_G g\tilde{G}g^{-1}.$$

Now $G$ acts on the Bass-Serre tree $T$ arising from this HNN
description. The stabilizers of the vertices are conjugates of
$\tilde G$ and the stabilizers of the edges are conjugates of
$K_1$ (and hence also of $K_2$). Let $u$ be a vertex of $T$ such
that $\tilde G=Stab(u)$, and let $e_1,e_2$ be two edges of $T$
incident at $u$ such that $K_1=Stab(e_1)$ and $K_2=Stab(e_2)$.

Now suppose that $g\notin\tilde{G}$.  Then $M_1\cap_G
gM_2g^{-1}\subset \tilde{G}\cap_G g\tilde{G}g^{-1}$ stabilises the
(nonempty) geodesic path $P$ in $T$ from $u$ to $g(u)$. Moreover,
since $\sigma(g,a_k)=0$, this path has even length and contains
the same number of forward-pointing and backward-pointing edges.
In particular, there is an intermediate vertex $v$ in $P$ which is
either the source of both the incident edges of $P$ or the target
of both the incident edges of $P$. We assume the latter.  (The
analysis of the former case is analogous.)

If $v=h(u)$, then the stabilisers of the edges of $P$ incident at
$v$ have the form $hsK_2s^{-1}h^{-1}$ and $htK_2t^{-1}h^{-1}$ for
some $s,t\in\tilde{G}$ with $s^{-1}t\notin K_2$.
By Bagherzadeh's Theorem \ref{b}, $sK_2s^{-1}\cap_{\tilde{G}}tK_2t^{-1}$
is cyclic, and hence the stabiliser of $P$ is cyclic, and the result
follows.

\medskip
Thus we are reduced to the case where $g\in\tilde{G}$.
But in that case $M_1$ and $\tilde{M_2}$ are Magnus subgroups of the one-relator group $\tilde{G}$, and
$$M_1\cap_G gM_2g^{-1} = M_1\cap_{\tilde{G}} g\tilde{M_2}g^{-1},$$
which is cyclic by Collins' Theorem \ref{col2}, unless
$$g\in M_1\cdot\tilde{M_2}\subseteq M_1\cdot M_2.$$

\end{proof}


\begin{thebibliography}{99}
\bibitem{Ba76}  G.H. Bagherzadeh, {\it Commutativity in one-relator groups}, J. London Math. Soc.(2),
{\bf 13}, 459--471(1976).
\bibitem{Bi81} S.D. Brodski\u\i, {\it Anomalous products of locally indicable groups} (Russian),
Algebraic systems,  51-77, Ivanov. Gos. Univ., Ivanovo, (1981).
\bibitem{Bog} O. Bogopolski, {\it A surface groups analogue of a theorem of Magnus}, in:
Geometric methods in group theory (ed J. Burillo, S. Cleary, M.
Elder, J. Taback and E. Ventura),  59--69, Contemp. Math., 372,
Amer. Math. Soc., Providence, RI, (2005).
\bibitem{BS} O. Bogopolski and K. Sviridov, A Magnus theorem for some
one-relator groups, preprint, (2006).
\bibitem{Col04} D.J. Collins, {\it Intersections of Magnus subgroups of one-relator groups},
Groups: topological, combinatorial and arithmetic aspects, L.
Math. Soc. Lec. Notes Ser., {\bf 311}, Cambridge University Press,
Cambridge, 255-296 (2004).
\bibitem{Col05} D.J. Collins, {\it Intersections of conjugates of Magnus subgroups of one-relator groups},
Zieschang Gedenkschrift (to appear).
\bibitem{Hem90} J. Hempel, {\it One-relator surface groups}, Math. Proc. Cambridge Philos. Soc.,
 {\bf 108}, 467-474 (1990).
\bibitem{How} J. Howie, {\it Some results on one-relator surface groups}, Bol. Soc. Mat. Mexicana
(3), {\bf 10}, Special Issue, 255-262 (2004); {\it Erratum}, {\it
ibid.}, 545--546 (2004).
\bibitem{How05} J. Howie, Magnus intersections in one-relator products, Michigan Math. Journal 53, 597--623 (2005).
\bibitem{Mg30}  W. Magnus, {\it \"{U}ber diskontinuierliche Gruppen mit einer definierenden Relation},
J. Reine Angew. Math., {\bf 163}, 141-165 (1930).
\bibitem{Mg31}  W. Magnus, {\it Untersuchungen \"{u}ber einige unendliche diskontinuierliche Gruppen}.
Math. Ann., {\bf 105}, 52-74 (1931).
\bibitem{New68} B.B. Newman, {\it Some results on one-relator groups}, Bull. Amer. Math. Soc., {\bf 74}, 568-571 (1968).
Math. Ann., {\bf 78}, 385-397 (1918).
\bibitem{MSS} M.S. Saeed, {\it One-relator quotients of surface groups}, PhD thesis, Heriot-Watt University, UK (2007).
\end{thebibliography}
\end{document}